\documentclass[twoside,12pt]{article}
\usepackage[]{amsmath,amsthm,amssymb}
\usepackage[latin1]{inputenc}

\begin{document}

\title{{\Large\bf  Discrete Lebedev's index  transforms}}

\author{Semyon  YAKUBOVICH}
\maketitle

\markboth{\rm \centerline{ Semyon   YAKUBOVICH}}{}
\markright{\rm \centerline{LEBEDEV'S  DISCRETE TRANSFORMS}}

\begin{abstract} {\noindent Discrete analogs of the Lebedev transforms with the product of the modified Bessel functions are introduced and investigated. Several expansions of suitable  functions and sequences in terms of the series and integrals, involving the modified and incomplete Bessel functions  are established. }

\end{abstract}
\vspace{4mm}



\vspace{4mm}

\section {Introduction and preliminary results}

In 1962 N.N. Lebedev    proved  (cf. [2]) the  following expansion  for all $x >0$

$$f(x) = - {4\over \pi^2} {d\over dx}  \int_0^\infty \tau \sinh(\pi\tau) K^2_{i\tau} (x)  \int_0^\infty K_{i\tau}(y) \left[I_{i\tau} (y) + I_{-i\tau} (y)\right] f(y) dy d\tau,\eqno(1.1)$$
where $I_\nu(z), K_\nu(z)$ are modified Bessel functions of the first and second kind, respectively, and $x^{-1/2} f(x)  \in L_1(0,1),\ x^{1/2} f(x) \in L_1(1,\infty)$.  It generates the reciprocal pair of the index transforms [5] with the product of the modified Bessel functions

$$g(\tau) =  \int_0^\infty K_{i\tau}(y) \left[I_{i\tau} (y) + I_{-i\tau} (y)\right] f(y) dy,\eqno(1.2)$$

$$f(x) = - {4\over \pi^2} {d\over dx}  \int_0^\infty \tau \sinh(\pi\tau) K^2_{i\tau} (x) g(\tau) d\tau.\eqno(1.3)$$
The main aim of the present paper is to introduce discrete analogs of the Lebedev transforms (1.2), (1.3) and study their mapping and inversion properties. Precisely, we will investigate the following transformations

$$a_n  =  {\pi\over 2\cosh(\pi n/2)} \int_0^\infty K_{in/2}(y) \left[I_{in/2} (y) + I_{-in/2} (y)\right] f(y) dy,\ n \in \mathbb{N},\eqno(1.4)$$

$$a_n  =  {1\over 2} \int_0^\infty K^2_{in/2}(y)  f(y) dy,\quad n \in \mathbb{N},\eqno(1.5)$$

$$f(x)  =  {\pi\over 2}  \sum_{n=1}^\infty {a_n\over \cosh(\pi n/2)} K_{in/2}(x) \left[I_{in/2} (x) + I_{-in/2} (x)\right], \ x >0,\eqno(1.6)$$

$$f(x)   =  {1\over 2} \sum_{n=1}^\infty a_n  K^2_{in/2}(x), \ x >0.\eqno(1.7)$$
In the sequel, inversion theorems will be proved for transformations (1.4)-(1.7), employing the theory of the discrete Kontorovich-Lebedev transform recently developed by the author (see [7]).  

As is known,  the modified Bessel function of the second kind or Macdonald function $K_\nu(z)$  is represented, for instance,  by the integral
(cf. [4], Vol. I,  Entry  2.4.18.4) 

$$K_{\nu}(z)= \int_{0}^\infty e^{-z\cosh (u) } \cosh(\nu u) du,\   {\rm Re}  z > 0, \  \nu \in \mathbb{C}.\eqno(1.8)$$
It satisfies the ordinary differential equation
$$  z^2{d^2u\over dz^2}  + z{du\over dz} - (z^2+\nu^2)u = 0,\eqno(1.9)$$
for which it is the solution that remains bounded as $z$ tends to infinity on the real line. It has the asymptotic behavior [6]
$$ K_\nu(z) = \left( \frac{\pi}{2z} \right)^{1/2} e^{-z} [1+ O(1/z)], \qquad z \to \infty,\eqno(1.10)$$
and near the origin
$$ K_\nu(z) = O\left ( z^{-|{\rm Re}\nu|}\right), \ z \to 0,\eqno(1.11)$$
$$K_0(z) = -\log z + O(1), \ z \to 0. \eqno(1.12)$$
Worth  mentioning is  the Lebedev inequality for the modified Bessel function (see  [6], p.219)

$$ \left| K_{i\tau} (x) \right| \le A \ {x^{-1/4} \over \sqrt{ \sinh (\pi\tau)} },\quad x, \tau > 0,\eqno(1.13)$$
where $A > 0$ is an absolute constant.   The kernels of discrete transforms (1.4)-(1.7) have representations in terms of the Erd{\' e}lyi-Kober integrals [6] (see [4], Vol. II, Entries 2.16.3.1, 2.16.3.6)

$${\pi\over 2\cosh(\pi n/2)}  K_{in/2}(x) \left[I_{in/2} (x) + I_{-in/2} (x)\right] = \int_0^x {K_{in} (t)\over (x^2-t^2)^{1/2}}\ dt,\eqno(1.14)$$

 $$  {1\over 2} \ K^2_{in/2}(x) =   \int_x^\infty  {K_{in} (t)\over (t^2-x^2)^{1/2}}\ dt.\eqno(1.15)$$
In the sequel we will also employ  the modified Struve functions ${\bf M}_\nu(z),\ {\bf L}_\nu(z)$, having the following relation with the modified Bessel function of the first kind  (cf. [3], Entry 11.2.6)

$$ {\bf M}_\nu(z) = \ {\bf L}_\nu(z) - I_\nu(z),\eqno(1.16)$$ 
and ${\bf L}_\nu(z) = - i e^{-\pi i\nu/2} {\bf H}_\nu(iz)$, where ${\bf H}_\nu(z)$ is the Struve function (see [3], Entry 11.2.2).

\section{Inversion theorems} 

We begin with

{\bf Theorem 1}.   {\it Let $f$ be a complex-valued function on $\mathbb{R}_+$ which is represented by the integral 

$$f(x) =  { 2 x \over \pi} \int_{-\pi}^\pi K_0(x \cosh(u)) \ \varphi(u) \cosh(u) du,\quad x >0,\eqno(2.1)$$ 
where $ \varphi(u) = \psi(u)\sinh(u)$ and $\psi$ is a  $2\pi$-periodic function, satisfying the Lipschitz condition on $[-\pi, \pi]$, i.e.

$$\left| \psi(u) - \psi(v)\right| \le C |u-v|, \quad  \forall \  u, v \in  [-\pi, \pi],\eqno(2.2)$$
where $C >0$ is an absolute constant.  Then the following inversion formula for  transformation $(1.4)$  holds

$$ f(x)  =  {1\over  \pi^3} \   \sum_{n=1}^\infty    \sinh(\pi n)  \Phi_n (x) a_n,\quad x > 0,\eqno(2.3)$$
where}

$$\Phi_n(x)= x  \int_{-\pi}^\pi K_0(x \cosh(u)) \ \sinh(2u) \sin(nu) du,\quad x >0,\ n \in \mathbb{N}.\eqno(2.4)$$

\begin{proof}   Plugging the right-hand side of the representation (1.14) in (1.4), we change the order of integration to obtain 

$$a_n  =  \int_0^\infty K_{in}(t ) \int_t^\infty  {f(y) \over (y^2-t^2)^{1/2}}\ dy dt.\eqno(2.5)$$
The justification of this interchange comes from the Fubini theorem via the use of H{\" o}lder's  inequality, generalized Minkowski's  inequality, asymptotic behavior (1.10)-(1.12) for the modified Bessel function and integral (2.1).   In fact, we derive

$$\int_0^\infty \left|K_{in}(t )\right|  \int_t^\infty  {|f(y) |\over (y^2-t^2)^{1/2}}\ dy dt\le \left( \int_0^\infty \left( \int_t^\infty  {|f(y) |\over (y^2-t^2)^{1/2}}\ dy \right)^p dt \right)^{1/p} $$

$$\times \left(  \int_0^\infty K^q_{0}(t ) dt\right)^{1/q}\le \int_0^\infty |f(y)| \ y^{-1/q} dy \left( \int_0^1  {1\over (1-t^2)^{p/2}}\ dt \right)^{1/p}  $$

$$\times \left(  \int_0^\infty K^q_{0}(t ) dt\right)^{1/q} \le \left( \int_0^1  {1\over (1-t^2)^{p/2}}\ dt \right)^{1/p} \left(  \int_0^\infty K^q_{0}(t ) dt\right)^{1/q}$$

$$\times   { 1 \over \pi}  \int_0^\infty y^{1/p} K_0( y ) dy  \int_{-\pi}^\pi  |\varphi(u)|  [ \cosh(u)]^{-1/p}  du $$

$$=  {1\over 2 \pi} B^{1/p}\left({1\over 2}, 1- {p\over 2}\right) \Gamma^2 \left({1\over 2}\left( 1 + {1\over p} \right) \right)  \left(  \int_0^\infty K^q_{0}(t ) dt\right)^{1/q}$$

$$\times   \int_{-\pi}^\pi  |\varphi(u)|  [ \cosh(u)]^{-1/p}  du < \infty,$$
where $1 <  p < 2,\ q = p/(p-1)$ and $B(a,b), \ \Gamma(z)$ are Euler's beta and gamma functions, respectively, [1].  For  the same reasons the inner integral with respect to $y$ in (2.5) can be written, using (2.1), Entry  3.14.1.9 in [1] and particular cases of the hypergeometric function ${}_0F_1$ (cf. [5], Vol. II),  in the form

$$\int_t^\infty  {f(y) \over (y^2-t^2)^{1/2}}\ dy =   {2\over \pi} \int_{-\pi}^\pi  \ \varphi(u) \cosh(u)  \int_t^\infty  { y  K_0(y \cosh(u)) \over (y^2-t^2)^{1/2}}\ dy  du $$

$$=   {2 \over \pi} \int_{-\pi}^\pi  \ \varphi(u) \cosh(u) \left[ - {\pi t\over 2} \  {}_0F_1\left( {3\over 2}; \  \left({t\cosh(u)\over 2}\right)^2 \right) \right.$$

$$\left. + {\pi\over 2\cosh(u)} \  {}_0F_1\left( {1\over 2}; \  \left({t\cosh(u)\over 2}\right)^2 \right) \right] du$$

$$=  {2 \over \pi} \int_{-\pi}^\pi  \ \varphi(u) \cosh(u) \left[ - {\pi \over 2\cosh(u) }  \sinh(t\cosh(u))  \right.$$

$$\left. + {\pi\over 2\cosh(u)} \cosh(t\cosh(u))  \right] du =   \int_{-\pi}^\pi e^{-t\cosh(u)}  \varphi(u) du.$$
Therefore, following the same scheme as in the proof of Theorem 5 in [7], we return to (2.5), substituting the latter expression,  changing the order of integration owing to the absolute and uniform convergence and  employ the formula (see [4], Vol. II, Entry 2.16.6.1)

$$ \int_0^\infty  e^{-x\cosh(u)} K_{in}(x)  dx = {\pi \sin( nu) \over \sinh (u) \sinh(\pi n)}\eqno(2.6)$$
to get finally

$$a_n  = {\pi\over  \sinh(\pi n)} \int_{-\pi}^\pi   \varphi(u) {\sin( nu) \over \sinh(u)} du.\eqno(2.7)$$
Let $S_N(x) $ denote  a partial sum of the series (2.3). Then, substituting  the value of $a_n$ by integral (2.7) and  $\Phi_n(x)$ by (2.4),  it gives 
$$S_N(x)  ={x\over \pi^2}  \sum_{n=1}^N   \int_{-\pi}^\pi K_0(x \cosh(t)) \ \sinh(2t) \sin(nt) dt \int_{-\pi}^\pi   {\varphi(u) \over \sinh(u)} \sin(nu) du.\eqno(2.8)$$
Hence,  calculating the sum via the known identity and invoking the definition of $\varphi$, equality (2.8) becomes

$$  S_N(x)  =  {x\over 4 \pi^2} \   \int_{-\pi}^\pi  K_0(x \cosh(t))  \sinh(2t)   \int_{-\pi}^{\pi}   {\varphi(u)+ \varphi(-u)  \over \sinh(u)} \  {\sin \left((2N+1) (u-t)/2 \right)\over \sin( (u-t) /2)}  du dt $$

$$=   {x\over 4 \pi^2} \   \int_{-\pi}^\pi   K_0(x \cosh(t)) \sinh(2t)   \int_{-\pi}^{\pi}  \left[ \psi(u)- \psi(-u) \right]  \  {\sin \left((2N+1) (u-t)/2 \right)\over \sin( (u-t) /2)}  du dt.\eqno(2.9)$$
Since $\psi$ is $2\pi$-periodic, we treat  the latter integral with respect to $u$ as follows 

$$  \int_{-\pi}^{\pi}  \left[ \psi(u)- \psi(-u) \right]  \  {\sin \left((2N+1) (u-t)/2 \right)\over \sin( (u-t) /2)}  du $$

$$=  \int_{ t-\pi}^{t+ \pi}  \left[ \psi(u)- \psi(-u) \right]  \  {\sin \left((2N+1) (u-t)/2 \right)\over \sin( (u-t) /2)}  du $$

$$=  \int_{ -\pi}^{\pi}  \left[ \psi(u+t)- \psi(-u-t) \right]  \  {\sin \left((2N+1) u/2 \right)\over \sin( u /2)}  du. $$
Moreover,

$$ {1\over 2\pi} \int_{ -\pi}^{\pi}  \left[ \psi(u+t)- \psi(-u-t) \right]  \  {\sin \left((2N+1) u/2 \right)\over \sin( u /2)}  du - \left[ \psi(t)- \psi(-t) \right] $$

$$=  {1\over 2\pi} \int_{ -\pi}^{\pi}  \left[ \psi(u+t)- \psi(t) + \psi (-t) - \psi(-u-t) \right]  \  {\sin \left((2N+1) u/2 \right)\over \sin( u /2)}  du.$$
When  $u+t > \pi$ or  $u+t < -\pi$ then we interpret  the value  $\psi(u+t)- \psi(t)$ by  formulas

$$\psi(u+t)- \psi(t) = \psi(u+t-2\pi)- \psi(t - 2\pi),$$ 

$$\psi(u+t)- \psi(t) = \psi(u+t+ 2\pi)- \psi(t +2\pi),$$ 
respectively.  Analogously, the value  $\psi(-u-t)- \psi(-t)$  can be treated.   Then   due to the Lipschitz condition (2.2) we have the uniform estimate
for any $t \in [-\pi,\pi]$

$${\left|  \psi(u+t)- \psi(t) + \psi (-t) - \psi(-u-t) \right| \over | \sin( u /2) |}  \le 2C \left| {u\over \sin( u /2)} \right|.$$
Therefore,  owing to the Riemann-Lebesgue lemma

$$\lim_{N\to \infty } {1\over 2\pi} \int_{ -\pi}^{\pi}  \left[ \psi(u+t)- \psi(-u-t)  - \psi(t) + \psi (-t) \right]  \  {\sin \left((2N+1) u/2 \right)\over \sin( u /2)}  du =  0\eqno(2.10)$$
for all $ t\in [-\pi,\pi].$    Besides, returning to (2.9), we estimate the iterated integral 

$$ \int_{-\pi}^\pi K_0(x \cosh(t)) | \sinh(2t) |  \int_{ -\pi}^{\pi} \left| \left[ \psi(u+t)- \psi(-u-t)  - \psi(t) + \psi (-t) \right]  \  {\sin \left((2N+1) u/2 \right)\over \sin( u /2)}  \right| du dt$$

$$\le  4 C \int_{0}^\pi K_0(x \cosh(t))  \sinh(2t)   dt   \int_{ -\pi}^{\pi}   \left| {u\over \sin( u /2)} \right| du < \infty,\ x >0.$$
Consequently, via  the dominated convergence theorem it is possible to pass to the limit when $N \to \infty$ under the  integral sign, and recalling (2.10), we derive

$$  \lim_{N \to \infty}   {x\over 4 \pi^2}  \int_{-\pi}^\pi K_0(x \cosh(t))  \sinh(2t)   \int_{ -\pi}^{\pi}  \left[ \psi(u+t)- \psi(-u-t)  - \psi(t) + \psi (-t) \right] $$

$$\times  \  {\sin \left((2N+1) u/2 \right)\over \sin( u /2)}  du dt =  {x\over 4 \pi^2}  \int_{-\pi}^\pi K_0(x \cosh(t))  \sinh(2t)   $$

$$ \times \lim_{N \to \infty}  \int_{ -\pi}^{\pi}  \left[ \psi(u+t)- \psi(-u-t)  - \psi(t) + \psi (-t) \right]  \  {\sin \left((2N+1) u/2 \right)\over \sin( u /2)}  du dt = 0.$$
Hence, combining with (2.9),  we obtain  by virtue of  the definition of $\varphi$ and $f$

$$ \lim_{N \to \infty}  S_N(x) =   {x\over \pi} \   \int_{-\pi}^\pi  K_0(x \cosh(t))  \cosh(u)   \left[ \varphi (t)+ \varphi(-t) \right] dt = f(x),$$
where the integral (2.1) converges since $\varphi \in C[-\pi,\pi]$.  Thus we established  (2.3), completing the proof of Theorem 1.
 
\end{proof} 

Concerning transformation (1.5), we have

{\bf Theorem 2}.   {\it Let $f$ be a complex-valued function on $\mathbb{R}_+$ which is represented by the integral 

$$f(x) =  \int_{-\pi}^\pi  \varphi(u)  \left[ {2\over \pi} +   x\cosh(u) {\bf M}_0\left(x\cosh(u)\right)\right] du,\quad x >0,\eqno(2.11)$$ 
where $ \varphi(u) = \psi(u)\sinh(u)$ and $\psi$ is a  $2\pi$-periodic function, satisfying the Lipschitz condition $(2.2)$.  Then the following inversion formula for  transformation $(1.5)$  takes place

$$ f(x)  =  {1\over  \pi^2} \   \sum_{n=1}^\infty    \sinh(\pi n)  \Psi_n (x) a_n,\quad x > 0,\eqno(2.12)$$
where}

$$\Psi_n(x)=    \int_{-\pi}^\pi \left[ {2\over \pi}  +   x\cosh(u) {\bf M}_0\left(x\cosh(u)\right)\right] \ \sinh(u) \sin(nu) du,\quad x >0,\ n \in \mathbb{N}.\eqno(2.13)$$

\begin{proof}  In the same manner as in the proof of Theorem 1 we substitute the right-hand side of (1.15) into (1.5), changing the order of integration. Hence we get 

$$a_n  =  \int_0^\infty K_{in}(t ) \int_0^t  {f(y) \over (t^2-y^2)^{1/2}}\ dy dt.\eqno(2.14)$$
Noting the integral representation for the modified Struve function (cf. [3], Entry  11.5.4)

$${\bf M}_0(z) = - {2\over \pi} \int_0^{\pi/2}  e^{-z\cos(\theta)} d\theta,$$
we easily find the estimate $| {\bf M}_0(z) | \le 1,\ {\rm Re} z \ge 0$. Hence the motivation of the interchange in (2.14) follows  immediately  from (2.11) and  Fubini's  theorem because

$$\int_0^\infty \left|K_{in}(t )\right|  \int_0^t  {|f(y) |\over (t^2-y^2)^{1/2}}\ dy dt\le  {2\over \pi} \int_{-\pi}^\pi | \varphi(u) | du  \int_0^\infty K_{0}(t ) dt   \int_0^1  {1\over (1-y^2)^{1/2}}\ dy $$

$$+   \int_{-\pi}^\pi | \varphi(u) | \cosh (u)  du  \int_0^\infty K_{0}(t ) t dt   \int_0^1  {y\over (1-y^2)^{1/2}}\ dy $$

$$=   \int_{-\pi}^\pi | \varphi(u) | \left[ {\pi\over 2} +  \cosh(u)\right] du < \infty.$$
Meanwhile, the inner integral with respect to $y$ in (2.14) can be represented via (1.16),  (2.11) and Entries 3.13.1.2,  3.15.1.4 in [1].  Therefore, invoking particular cases of the hypergeometric functions ${}_0F_1,\  {}_1F_2$  (see [3], Vol. III, Entries 7.13.1.6,  7.14.2.78) we find 

$$ \int_0^t  {f(y) \over (t^2-y^2)^{1/2}}\ dy =  \int_{-\pi}^\pi  \varphi(u) \left[ 1 -   \cosh(u)   \int_0^t  {I_0(y \cosh(u)) y \over (t^2-y^2)^{1/2}}\ dy \right.$$

$$\left. +  \cosh(u)   \int_0^t  {{\bf L}_0(y \cosh(u)) y \over (t^2-y^2)^{1/2}}\ dy\right] du =  \int_{-\pi}^\pi  \varphi(u) \left[ 1 -   t \cosh(u)\  {}_0F_1\left({3\over 2} ; \  \left[{t\cosh(u)\over 2}\right]^2 \right) \right.$$

$$\left. +  {[t \cosh(u)]^2\over 2}  \     {}_1F_2\left( 1;\ {3\over 2},\ 2  ; \  \left[ {t\cosh(u)\over 2}\right]^2 \right)  \right] du  =   \int_{-\pi}^\pi  \varphi(u) \left[ \cosh ( t\cosh(u))\right.$$ 

$$\left. -  \sinh (t \cosh(u)) \right] du =  \int_{-\pi}^\pi e^{-t\cosh(u) } \varphi(u) du.\eqno(2.15)$$
Thus, returning to (2.14) and appealing to (2.6),  the discrete transformation (1.5) takes the form (2.7).  Now in the same manner as in the proof of Theorem 1 we derive inversion formula (2.12) with the kernel (2.13).

\end{proof} 

The following result  establishes the inversion formula for the discrete transformation (1.6).   Indeed, we have

{\bf Theorem 3}. {\it   Let the sequence $\{a_n\}_{n\in \mathbb{N}}$ be such that the following series converges

$$\sum_{n=1}^\infty  |a_n|  e^{-\pi n/2} < \infty.\eqno(2.16)$$
Then the discrete transformation $(1.6)$ can be inverted by the formula

$$a_n =  { 1 \over \pi^3} \ \sinh(\pi n) \int_0^\infty   \Phi_n(x) f(x) dx,\ n \in \mathbb{N},\eqno(2.17)$$
where the kernel $\Phi_n(x)$ is defined by $(2.4)$ and  integral  $(2.17)$ converges absolutely.}

\begin{proof} In fact,  substituting (1.6), (1.14) and (2.4) on the right-hand side of (2.17), we estimate the integral, recalling inequality (1.13), to obtain 

$$ \int_0^\infty   \left|\Phi_n(x) f(x)\right| dx = \int_0^\infty x \left|  \int_{-\pi}^\pi K_0(x \cosh(u)) \ \sinh(2u) \sin(nu) du \right.$$

$$\left. \times \sum_{m=1}^\infty a_m \int_0^x {K_{im} (t)\over (x^2-t^2)^{1/2}}\  dt \right| dx \le  4 A \int_0^\infty x^{3/4}  \int_{0}^\pi K_0(x \cosh(u)) \ \sinh(2u) du dx $$

$$ \times \sum_{m=1}^\infty |a_m| e^{-\pi m/2}  \int_0^1 {t^{-1/4} \over (1-t^2)^{1/2}} dt =  2^{15/4} A \sqrt\pi \  \Gamma \left({3\over 8}\right) \Gamma \left({7\over 8}\right)$$

$$\times  \left(\cosh^{1/4}(\pi)-1\right)  \sum_{m=1}^\infty |a_m| e^{-\pi m/2} < \infty,$$
owing to assumption (2.16).  Therefore the interchange of the order of integration and summation is allowed to get

$$  { 1 \over \pi^3} \ \sinh(\pi n) \int_0^\infty   \Phi_n(x) f(x) dx =    { 1 \over \pi^3} \ \sinh(\pi n)  \sum_{m=1}^\infty a_m \int_{-\pi}^\pi  \ \sinh(2u) \sin(nu)$$

$$   \int_0^\infty K_{im} (t)  \int_t^\infty { x K_0(x \cosh(u)) \over  (x^2-t^2)^{1/2}} dx dt du.\eqno(2.18) $$
The integral with respect to $x$ in (2.18) is calculated in the proof of Theorem 1, and we have the formula 

$$  \int_t^\infty { x K_0(x \cosh(u)) \over  (x^2-t^2)^{1/2}} dx = {\pi  e^{-t\cosh(u)}\over 2\cosh(u)},\quad t > 0.$$
Hence, recalling (2.6), we  derive finally from (2.18)

$$   { 1 \over \pi^3} \ \sinh(\pi n) \int_0^\infty   \Phi_n(x) f(x) dx =    { 1 \over \pi} \ \sinh(\pi n)  \sum_{m=1}^\infty {a_m  \over \sinh(\pi m)} $$

$$\times \int_{-\pi}^\pi   \sin(nu) \sin(m u) du = a_n.$$

\end{proof}

Finally, we demonstrate the inversion theorem for the discrete transform (1.7). 

{\bf Theorem 4}. {\it   Let the sequence $\{a_n\}_{n\in \mathbb{N}} \in l_1$, i.e.  the  series 

$$\sum_{n=1}^\infty  |a_n|  < \infty\eqno(2.19)$$
converges. Then for the discrete transformation $(1.7)$ the following inversion formula holds 

$$a_n =  { 1 \over \pi^2} \ \sinh(\pi n) \int_0^\infty   \Psi_n(x) f(x) dx,\ n \in \mathbb{N},\eqno(2.20)$$
where the kernel $\Psi_n(x)$ is defined by $(2.13)$ and  integral  $(2.20)$ converges absolutely.}

\begin{proof} Indeed, similarly to the proof of Theorem 3 we recall (1.7), (1.15), (2.13) to find the estimate

$$\int_0^\infty  \left|  \Psi_n(x) f(x)\right| dx =   \int_0^\infty  \left| \int_{-\pi}^\pi \left[ {2\over \pi}  +   x\cosh(u) {\bf M}_0\left(x\cosh(u)\right)\right] \ \sinh(u) \sin(nu) du\right.$$

$$\left.  \sum_{m=1}^\infty a_m  \int_x^\infty  {K_{im} (t)\over (t^2-x^2)^{1/2}}\ dt \right| dx \le  2 \int_0^\infty   \int_{0}^\pi \left[ {2\over \pi}  +   x\cosh(u)\left| {\bf M}_0\left(x\cosh(u)\right)\right|\right] \ \sinh(u) du$$

$$\times  \sum_{m=1}^\infty  |a_m|  \int_1^\infty  {K_{0} (xt)\over (t^2- 1)^{1/2}}\ dt dx =    2 \sum_{m=1}^\infty  |a_m|   \int_1^\infty  {1\over t (t^2- 1)^{1/2}} \int_0^\infty  K_{0} (x) \int_{0}^\pi \left[ {2\over \pi}  +   {x\over t}\cosh(u) \right.$$

$$\left. \left| {\bf M}_0\left({x\cosh(u)\over t} \right)\right|\right] \ \sinh(u) du dx dt  \le  2\sum_{m=1}^\infty  |a_m|   \int_1^\infty  {1\over t (t^2- 1)^{1/2}} \int_0^\infty  K_{0} (x) \int_{0}^\pi \left[ {2\over \pi}  +   {x\over t}\cosh(u) \right]$$

$$\times  \ \sinh(u) du dx dt  = \left[ 2 \pi \sinh^2\left({\pi\over 2}\right) +  \sinh^2(\pi) \right]  \sum_{m=1}^\infty  |a_m| < \infty.$$
Therefore, taking into account (2.15),  it says

$$ \int_0^t  { 2/ \pi  +   x\cosh(u) {\bf M}_0\left(x\cosh(u)\right) \over (t^2-x^2)^{1/2}}\ dx = e^{- t\cosh(u)},\quad t > 0. $$
Thus we  derive via (2.6)

$$ { 1 \over \pi^2} \ \sinh(\pi n) \int_0^\infty   \Psi_n(x) f(x) dx =   { 1 \over \pi^2} \ \sinh(\pi n) \sum_{m=1}^\infty a_m   \int_{-\pi}^\pi  \sinh(u) \sin(nu)   \int_0^\infty   K_{im} (t) $$

$$\times  \int_0^t   \left[ {2\over \pi}  +   x\cosh(u) {\bf M}_0\left(x\cosh(u)\right)\right]   { dx dt du \over (t^2-x^2)^{1/2}} $$

$$=   { 2 \over \pi} \ \sinh(\pi n) \sum_{m=1}^\infty {a_m \over \sinh(\pi m)}    \int_{0}^\pi   \sin(nu) \sin(mu) du = a_n.$$

\end{proof}

\bigskip
\centerline{{\bf Acknowledgments}}
\bigskip

\noindent The work was partially supported by CMUP, which is financed by national funds through FCT (Portugal)  under the project with reference UIDB/00144/2020.

\bigskip
\centerline{{\bf References}}
\bigskip
\baselineskip=12pt
\medskip
\begin{enumerate}

\item[{\bf 1.}\ ]  Yu.A. Brychkov, O.I. Marichev, N.V. Savischenko,   {\it Handbook of Mellin Transforms.}  Advances in Applied Mathematics, CRC Press,  Boca Raton, 2018.

\item[{\bf 2.}\ ] N.N. Lebedev,  On an integral representation of an arbitrary function in terms of squares of Macdonald functions with imaginary index. {\it Sibirsk. Mat. Zh.} {\bf 3} ( 1962),  213- 222 (in Russian). 

\item[{\bf 3.}\ ] NIST Digital Library of Mathematical Functions. http://dlmf.nist.gov/, Release 1.0.17 of 2017-12-22. F. W. J. Olver, A. B. Olde Daalhuis, D. W. Lozier, B. I. Schneider, R. F. Boisvert, C. W. Clark, B. R. Miller and B. V. Saunders, eds.

\item[{\bf 4.}\ ] A.P. Prudnikov, Yu.A. Brychkov and O.I. Marichev, {\it Integrals and Series}. Vol. I: {\it Elementary
Functions}, Vol. II: {\it Special Functions}, Gordon and Breach, New York and London, 1986, Vol. III : {\it More special functions},  Gordon and Breach, New York and London,  1990.

\item[{\bf 5.}\ ] S.G. Samko,  A. A. Kilbas and O.I. Marichev,  {\it Fractional integrals and derivatives. Theory and applications. }  Gordon and Breach Science Publishers, Yverdon, 1993. 

\item[{\bf 6.}\ ] S. Yakubovich, {\it Index Transforms}, World Scientific Publishing Company, Singapore, New Jersey, London and
Hong Kong, 1996.

\item[{\bf 7.}\ ]  S. Yakubovich, Discrete Kontorovich-Lebedev transforms, ArXiv: 1908.01392.

\end{enumerate}

\vspace{5mm}

\noindent S.Yakubovich\\
Department of  Mathematics,\\
Faculty of Sciences,\\
University of Porto,\\
Campo Alegre st., 687\\
4169-007 Porto\\
Portugal\\
E-Mail: syakubov@fc.up.pt\\

\end{document}